\documentclass[12pt,a4paper]{article}
 \usepackage{geometry,bbm}
\usepackage{fullpage}
\usepackage{mathrsfs}
\usepackage{dsfont}
\usepackage[utf8]{inputenc}
\usepackage{amssymb}
\usepackage{amsmath}
\usepackage{amsthm}
\usepackage[round]{natbib}
\usepackage{url}
\usepackage{graphicx}
\usepackage{color}
\usepackage[colorlinks=true,linkcolor=blue,citecolor=blue,pdfborder={0 0 0}]{hyperref}
\usepackage{bm}

\usepackage{paralist}
\theoremstyle{plain}
\newtheorem{prop}{Proposition}[section]

\newtheorem{theorem}[prop]{Theorem}

\newtheorem{lemma}[prop]{Lemma}

\newtheorem{corollary}[prop]{Corollay}

\theoremstyle{remark}
\newtheorem{remark}[prop]{Remark}

\numberwithin{equation}{section}

\title{On the strong approximation and functional limit laws for the increments of the non-overlapping
$k$-spacings processes}

\author{Salim Bouzebda\footnote{Laboratoire de Math\'ematiques Appliqu\'ees, Universit\'e de Technologie de Compi\`egne, B.P. 529, 60205 Compi\`egne cedex, France, {E-mails:} \texttt{salim.bouzebda@utc.fr, nabil.nessigha@gmail.fr}} \,  and Nabil Nessigha}

\begin{document}
\maketitle

\begin{abstract}
\noindent The first aim of the present paper, is to
establish  strong approximations of the uniform non-overlapping
$k$-spacings process extending the results of
\cite{alybeirlanthorvart1984}. Our methods rely on the invariance
principle in \cite{masonvanzwet1987}.
The second goal, is to generalize the \cite{Dindar1997} results for
the increments of the spacings quantile process to the uniforme
non-overlapping $k$-spacings quantile process.  We apply the last result to characterize the limit laws of functionals of the  increments $k$-spacings quantile process.
\medskip


\medskip

\end{abstract}

\section{Introduction}
\noindent Let $U_1,U_2,\ldots$,  be
independent and identically distributed $(i.i.d\cdot)$ uniform $[0,1]$
random variables $(r.v.^{,}s)$ defined on the same probability
space $(\Omega,A,\mathbb{P})$. Denote by $0=:U_{0,n}\leq U_{1,n}\leq \cdots
\leq U_{n-1,n}\leq U_{n,n}:= 1$, the order statistics of
$U_1,U_2,\ldots,U_{n-1}$, and $0,1$. The
corresponding non-overlapping {$k$}-spacings are then defined by
\begin{equation}
\begin{array}{rll}\label{www2}
D^{k}_{i,n} & := & U_{im,n}-U_{(i-1)m,n},~~\mbox{for}~~ 1 \leq i \leq N-1,\\
D^{k}_{N,n} & := & 1-U_{(N-1)m,n},\\
\end{array}
\end{equation}
\noindent where $ N=\lfloor n/k \rfloor$, with $\lfloor u \rfloor
\leq u < \lfloor u \rfloor +1$ denoting the  integer part of $u$.
When $k=1$, i.e., $N=n$, the $k$-spacings reduce to  the usual
$1$-spacings (or simple spacings) defined by
$$D_{i,n}^{1}=U_{i,n}-U_{i-1,n}, ~~\mbox{for}~~i=1,\ldots,n.$$Simple spacings
have received a great deal of attention in the literature. We
refer to  \cite{Deheuvels1986}, \cite{pyk1965,Pyke1972},
\cite{Shorack1972}, \cite{raosethurman1975}, \cite{Beirlant84}
and \cite{DeheuvelsMason1991} for
details. Throughout the sequel, $k\geq 1$ will
denote a fixed integer. In applications it is more convenient to
use the normalized non-overlapping $k$-spacings
$\{{k}D_{i,n}^{k}: 1\leq i \leq N\}$. For a fixed $k\geq 1$, as
$n\rightarrow \infty$, the distribution function of
${k}D_{i,n}^{k}$ (which is independent of the index $i$ with
$1\leq i \leq N-1$) converges to the distribution function
$F_{k}(\cdot)$, of a standard gamma random variable with expectation
{$k$}, given by
\begin{equation}\label{zzz1}
 F_{k}(t):=\frac{1}{(k-1)!}\int_0^t x^{k-1}e^{-x}dx=
\int_0^t F_{k}(t)dt, ~~\mbox{for}~~t\geq 0,
\end{equation}
\noindent where
\begin{equation}\label{zzz2}
F_{k}(t)= \frac{t^{k-1}e^{-t}}{(k-1)!}
~~\mbox{and}~~F_{k}(t)=0,~~\mbox{for}~~ t<0.
\end{equation}
For each choice of $k\geq 1$, the empirical $k$-spacings process
is defined by
\begin{equation}\label{emprocess}
\alpha_n(x):=N^{1/2}\left(\widehat{F}_n(x)-F_{k}(x)\right),~~\mbox{for}~~ x>0,
\end{equation}
where $\widehat{F}_n(\cdot)$ is the empirical distribution function of
$\{{k}D_{i,n}^{k}: 1\leq i \leq N\}$, defined for $n\geq m$, by
\begin{equation}
\widehat{F}_n(x):=\frac{1}{N}\sum_{i=1}^N
\mathds{1}_{\left\{{k}D_{i,n}^{k}\leq x\right\}},~~\mbox{for}~~ x\in
\mathbb{R},
\end{equation}
\noindent with $\mathds{1}_{\{A\}}$ denoting the indicator function of
the event A. We will need the following
additional notation and definitions. Let
$
M_{1:n}^{k}\leq M_{2:n}^{k}\leq \cdots \leq M_{N:n}^{k},
$
be the order statistics of $\{D_{i,n}^{k}: 1\leq i \leq N\}$.
The quantile $k$-spacings function is given by
$$\widehat{Q}_n(t):= \left\{
\begin{array}{lll}
{k}M_{i,n}^{k}, &\mbox{if} &\ \frac{i-1}{N}< t\leq \frac{i}{N},~~
 i=1,2,\ldots,N,
\\
0,  &\mbox{if}&\  t=0.
\end{array}
\right.$$

\noindent Let us introduce 
\begin{equation}
Q_{k}(t)=\inf\left\{x\geq 0: F_{k}(x)\geq t\right\},
\end{equation}
and $$f_{k}(t)=\frac{d}{dt}F_{k}(t).$$The quantile $k$-spacings
process $\{\gamma_n(t):0\leq t \leq 1\}$ is then defined by
\begin{equation}
\gamma_n(t):=N^{1/2}f_{k}\left(Q_{k}(t)\right)\left(Q_{k}(t)-\widehat{Q}_n(t)\right), ~~\mbox{for}~~
0\leq t \leq 1.
\end{equation}
 In \cite{Deheuvels1985}, more than 60 references are given on this subject, with statistical applications such as testing uniformity or goodness-of-fit tests. Weak convergence results for the process $\{\alpha_n(x): 0\leq x<\infty, n\geq 1 \}$ were obtained by \cite{pyk1965}, \cite{Shorack1972}, \cite{raosethurman1975}, \cite{Aly1983} and \cite{Beirlant84}. Further, \cite{Durbin1975} obtained tables for the limiting distribution of the Kolmogorov-Smirnov (K-S)  statistic based on $\{\alpha_n(x): 0\leq x<\infty, n\geq 1 \}$. Here, we mention that 
\cite{pyk1965} was the first to suggest the use of the K-S and Cram\'er-von Mises  functionals of
$\{\alpha_n(x): 0\leq x<\infty, n\geq 1 \}$.  \cite{delpinoguido1979} consider the $k$-spacings as in (\ref{www2}) and characterized the limiting distribution of the statistics 
$$
W_n(g,k)=N^{-1/2}\sum_{i=1}^N(g(NkD_{i,N}^k)-a),
$$
where $g(\cdot)$ is a smooth function, $k$ is fixed and $a=\mathbb{E}[g(Y)]$, $Y$ is a \textit{r.v.} with a density function $f_k(y)$. These statistics $W_n(g,k)$ can be used for testing goodness-of-fit to a uniform distribution.
For application of the spacing in statistical tests and others we may refer to \cite{Magnus2013, Magnus2008}, \cite{Tung2012}, \cite{Tung2012B}, \cite{Deheuvels2006A} and \cite{Baryshnikov2009}.
 It is worth noticing that  simple spacings have received a great deal of
attention in the literature, we may refer to  \cite{Deheuvels2011} and \cite{Alvarez2017} for details. In \cite{Deheuvels2011}, the author obtain an explicit description of the limiting Gaussian process generated by the sample spacings from a non-uniform distribution.  \cite{alybeirlanthorvart1984} obtained strong approximations results for the  empirical process and the quantile process based on non-overlapping $k$-spacings and also the weak convergence of these processes in $\Vert/q\Vert$-metrics, refer to the last reference for definition.

\noindent The aim in this paper is to obtain a refinement of the
strong approximation results for $\{\alpha_n(x): 0\leq x<\infty, n\geq 1 \}$ and $\{\gamma_n(t):0\leq t\leq1,n\geq1\}  $
obtained by \cite{alybeirlanthorvart1984}. Their main tool is the
well known (KMT) invariance principle introduced in \cite{KMT1975}. In our approach we shall
make use the refinement of the KMT inequality for the Brownian
bridge approximation of uniform empirical and quantile processes
presented respectively in \cite{masonvanzwet1987}. This approach is based on the approximation of
the ${k}$-spacings process in the interval  $[0,a]$, with $a\leq 1$. In order to prove the invariance principle, we use the
same method developed in \cite{alybeirlanthorvart1984}, which is
based on the following representation of simple spacings given by
\cite{pyk1965}.  In the sequel of this section,  we use  a notation similar   to that used  in \cite{alybeirlanthorvart1984}
including some changes absolutely necessary for our setting. Let $E_1,E_2,\ldots$ denote an
$i.i.d.$ sequence of exponential $r.v.^{,}s$ with mean $1$ and set
$$S_n:=\sum_{i=1}^nE_i.$$
 Then, for each $n>1$, we have the
distributional identity
\begin{equation}\label{spacingsidentete}
\left\{U_{i,n}-U_{i-1,n}:  1\leq i \leq n
\right\}{\stackrel{d}{=}}\left\{\frac{E_i}{S_n}: 1\leq i \leq
n\right\}.
\end{equation}
where $ {\stackrel{d}{=}} $ denote the distributional equality.
\noindent Consequently we obtain the following representation of
the non-overlapping $k$-spacings
\begin{eqnarray}\label{m-sacingsiden}
\lefteqn{\left\{D_{i,n}^{k}, 1\leq i \leq
N-1,D_{N,n}^{k}\right\}}\nonumber \\&
\stackrel{d}{=}&\left\{\left(\sum_{\ell=i}^{i+k-1}E_{\ell}\right)\big/
S_n,\i=1,k+1,\ldots,\left(\left\lfloor\frac{n}{k}
\right\rfloor-1\right)k+1,
\left(\sum_{\ell=k\lfloor\frac{n}{k}\rfloor+1}^{n}E_{\ell}\right)\big/S_n\right\}.
\end{eqnarray}
 In particular,
if $n=Nk$ is an integer multiple of ${k}$, then
\begin{equation}\label{m-sacingsiden(n={k})}
\left\{D_{i,n}^{k}, 1\leq i \leq N
\right\}{\stackrel{d}{=}}\left\{Y_i/T_N, 1\leq i \leq N \right\},
\end{equation}
\noindent where
\begin{equation}\label{Y_i}
Y_i:=\sum_{\ell=(i-1)k+1}^{ik}E_{\ell}, ~\mbox{for}~i=1,2,\ldots,N,
\end{equation}
is a sequence of $i.i.d.$ $r.v.^{,}s$
with distribution function $F_{k}(\cdot)$ and
$$T_N=\sum_{i=1}^NY_i.$$ Now, we denote by $G_N(\cdot)$
the empirical distribution function and by $K_N(\cdot)$ the empirical
quantile function of the sequence $Y_1,\ldots,Y_N$, respectively,
defined by
\begin{equation}
G_N(x):=\frac{1}{N}\sum_{i=1}^N\mathds{1}_{\{Y_i\leq x\}},
~\mbox{for}~x\in \mathbb{R}^+,
\end{equation}
and
\begin{equation}\label{K_N(t)}
K_N(t):=\inf\{x:G_N(x)\geq t\},~\mbox{for}~ 0\leq t \leq 1.
\end{equation}
\noindent Let $\{\beta_N(x): 0\leq x<\infty, N\geq 1\}$ and $\{\kappa_N(t):1\leq t \leq 1, N\geq  1\}$ be the corresponding
empirical and quantile processes, respectively, defined by
\begin{equation}
\beta_N(x):=\sqrt{N}\left(G_N(x)-F_{k}(x)\right), ~\mbox{for
 }~x\in \mathbb{R}^+,
\end{equation}
and
\begin{equation}
\kappa_N(t):=\sqrt{N}f_{k}\left(Q_{k}(t)\right)\left(Q_{k}(t)-K_N(t)\right),
~\mbox{for }~ 0\leq t \leq 1.
\end{equation}
By (\ref{m-sacingsiden(n={k})}), we have the following
representation
\begin{equation}\label{alpha}
\left\{\alpha_{{Nk}}(x), 0\leq x<\infty
\right\}{\stackrel{d}{=}}\left\{\alpha_N^1(x)=\beta_N\left(x\frac{T_N}{{Nk}}\right)+
\mathcal{R}_N(x), 0\leq x <\infty\right\},
\end{equation}
\noindent where
$$\mathcal{R}_N(x)=N^{1/2}\left(F_{k}\left(x\frac{T_N}{{Nk}}\right)
-F_{k}(x)\right).$$ \noindent In the same way, by
(\ref{m-sacingsiden(n={k})}), and the definition of the empirical
quantile function $K_N(\cdot)$, we have the following representation for the process 
$\left\{\gamma_{{Nk}}(t): 0\leq t <
1, N\geq 1\right\}$
\begin{equation}\label{gamma}
\left\{\gamma_{{Nk}}(t): 0\leq t <
1\right\}{\stackrel{d}{=}}\left\{\gamma_N^1(t)=\frac{{Nk}}{T_N}\left(
\kappa_N(t)+N^{1/2} \left(\frac{T_N}{{Nk}}-1\right)\phi_k(t)\right):
0\leq t < 1\right\},
\end{equation}
where
\begin{equation}
\phi_k(t)=f_{k}(Q_{k}(t))Q_{k}(t)\nonumber.
\end{equation}
The methodology  in the proof of the results of \cite{alybeirlanthorvart1984} will be essential to obtain our main theorems.  
In the present work, in the derivation of our asymptotic results we give a rigorous proof  with all mathematical details and give some constants explicitly.

The rest of the paper is organized as follows.  We establish  local strong approximations of the uniform non-overlapping
$k$-spacings process in the forthcoming section.  In Section \ref{sectionf}, we establish the functional limit laws for the increments of the quantile process of
non-overlapping $k$-spacings processes.  We provide an application concerning the limit laws of functionals of the $k$-spacings process. To avoid interrupting the flow of the
presentation, all mathematical developments are relegated to Section \ref{proofc}.
\section{Local Strong Approximation }
\subsection{Preliminaries}
Let us begin by introducing some Gaussian
processes playing a central role in strong approximations
theory. Let $W=
\{W(s): s\geq 0\}$ and $B=\{B(u) : u\in[0,1]\}$ be
the standard Wiener process and Brownian bridge, that is, the centered Gaussian
processes with continuous sample paths, and covariance functions
$$
\mathbb{E}(W(s)W(t))=s\wedge t,\quad\text{for}\quad s,t\geq 0
$$
and
$$
\mathbb{E}(B(u)B(v))=u\wedge v- uv,\quad\text{for}\quad u,v\in[0,1].
$$
The interested reader may refer to \cite{csorgorevesz1981} for details on
the Gaussian processes mentioned above.
In the sequel, the underlying probability space $(\Omega,A,\mathbb{P})$ is assumed
to be rich enough, in the sense that an independent sequence of
Gaussian processes, which is independent of the originally given
i.i.d. sequence of random vectors, can be constructed on this
probability space. This is a technical requirement which
allows for the construction of the Gaussian processes in our following Theorems. Since one can expand the
underlying probability space, this assumption is not restrictive,
refer to \cite{deAcosta1982}, \cite[Lemma A1]{Berkes1979} and  \cite{KMT1975} for further details. Throughout the paper we set $\log_{+}(u)=\log(u\vee e) , ~\mbox{for}~
u\in \mathbb{R}$. Let us recall the following theorem which is a weak version of the result in 
\cite{masonvanzwet1987}.\vskip5pt
\begin{theorem}\label{Masonvanzwet theorem}
(\cite{masonvanzwet1987}). There exists a sequence of empirical
processes $\beta_N$ based on $Y_1,\ldots,Y_N$ and a sequence of
Brownian bridges $\{B^{(1)}_N(t) :0\leq t \leq 1\} $ such that, for
all $\varepsilon
>0$ and $0 \leq a \leq\ 1$, we have
\begin{equation}
\mathbb{P}\left(\sup_{0 \leq x \leq
Q_{k}(a)}|\beta_N(x)-B^{(1)}_N(F_{k}(x))| \geq
\mathcal{A}N^{-1/2}(\log aN) \right)\leq
\mathcal{B}N^{-\varepsilon},
\end{equation}
where $\mathcal{A}$ and $\mathcal{B}$ are positive constants
depending on $\varepsilon$ and $a$.
\end{theorem}
A similar result is needed for the quantile
process $\{\kappa_n(t): 0\leq t\leq 1, n\geq 1\}$. For this, we consider deviations between the
quantile process $\{\kappa_N(t):0\leq t\leq1, N\geq 1 \}$ and the approximating Brownian bridges
$\{B^{(1)}_N(t),0\leq t \leq 1\} $ on $[0,a]$, instead of $[0,1]$.
We formulate this idea in the following theorem.
\begin{theorem}\label{1er theorem}
Let $\{B^{(1)}_N(t),0\leq t \leq 1\} $ be as in of Theorem
\ref{Masonvanzwet theorem}. Then for all $\varepsilon
>0$ and $n\geq m$, we have
\begin{equation}
\mathbb{P}\left(\sup_{0 \leq t \leq a}|\kappa_N(t)-B^{(1)}_N(t)| \geq
A_1N^{-1/4}(\log aN)^{3/4} \right)\leq B_1N^{-\varepsilon},
\end{equation}
 for all  $0 \leq a \leq\ 1$, where
 $A_{1}$ and  $ B_1$ are positive constants.
\end{theorem}
\noindent The proof of Theorem \ref{1er theorem} is postponed until Section \ref{proofc}.
\subsection{Main results}
\noindent  Our next theorem describes the strong approximations of the process $\left\{\gamma_{{Nk}}(t): 0\leq t <
1\right\}$.
\begin{theorem}\label{2eme theorem}
There exists a sequence $\left\{W_{{Nk}}(t): 0\leq t\leq
1, \geq 1\right\}$ of Gaussian processes, such that the following
properties hold. We have
\begin{equation}
\mathbb{E}W_{{Nk}}(t)=0,\nonumber
 \mbox{~and~}\mathbb{E}W_{{Nk}}(t)W_{{Nk}}(s)=\min(t,s)-ts-\frac{1}{n}\phi_{k}(t)\phi_{k}(s),\nonumber
\end{equation}
and
\begin{equation}
\phi_{k}(t)=f_{k}\left(Q_{k}(t)\right)Q_{k}(t)\nonumber.
\end{equation}
Moreover, for each $\varepsilon>0$, there exists constants
$A_{2}>0$ and $ B_{2}>0$, such that, for all $n\geq k$ and $a\in
[0,1]$ we have
\begin{equation}
\mathbb{P}\left(\sup_{0 \leq t \leq a }|\gamma_{{Nk}}(t)-W_{{Nk}}(t)| >
A_{2}N^{-1/4}(\log aN)^{3/4} \right)\leq
B_{2}N^{-\varepsilon}\nonumber.
\end{equation}
\end{theorem}
\noindent The proof of Theorem \ref{2eme theorem} is postponed until Section \ref{proofc}.\\
\noindent Our next theorem describes the strong approximations of the process $\{\alpha_n(x): 0\leq x <
\infty\}$.
\begin{theorem}\label{3eme theorem}
There exists a sequence of Gaussian processes $\{V_n(x): 0\leq
x\leq \infty, n\geq 1\}$, such that the following properties hold. We have
\begin{equation}
\mathbb{E}V_n(x)=0,
\end{equation}
and
\begin{equation}
\mathbb{E}V_n(x)V_n(y)=\min\left(F_{k}(x),F_{k}(y)\right)
-F_{k}(x)F_{k}(y)-\frac{1}{k}xyF_{k}(x)F_{k}(y).
\end{equation}
Moreover, for all $\varepsilon>0$ and $a\in [0,1]$ we have
\begin{equation}
\mathbb{P}\left(\sup_{0 \leq x \leq Q_{k}(a)}|\alpha_n(x)-V_n(x)| \geq
A_{3}N^{-1/4}(\log aN) \right)\leq B_{3}N^{-\varepsilon}\nonumber,
\end{equation}
where $A_{3}>0$ and $B_{3}>0$ are positive constants.
\end{theorem}
\noindent The proof of Theorem \ref{3eme theorem} is postponed until Section \ref{proofc}.\\
Immediate consequences of Theorems~\ref{2eme theorem} and \ref{3eme theorem} are  upper bounds for the
convergence of distributions of smooth functionals of $\{\gamma_{{Nk}}(t):0\leq t\leq 1, N\geq 1\}$ and $\{\alpha_n(x): 0\leq x< \infty, n\geq 1 \}$.
 Notice that the following corollary
is the analogous of the Corollary of \cite{KMT1975} page 113.
Let $\mathcal{D}(\mathbb{A})$ be the space of right-continuous real-valued
functions defined on $\mathbb{A}$ which have left-hand limits, equipped
with the Skorohod topology; refer to \cite{Billingsley1968} for further details.
\begin{corollary}\label{corol22}
Let $\Phi(\cdot)$ be a  functional defined
on  the space $\mathcal{D}(\mathbb{R})$, satisfying a Lipschitz condition
$$
|\Phi(v)-\Phi(w)|\leq L \sup_{t\in\mathbb{R}} |v(t)-w(t)|.
$$
Assume further that the distribution of the r.v. $\Phi\big(W_{{Nk}}(F(\cdot))\big)$
 has a bounded density. Then, as $n\to\infty$,
\begin{equation}\label{estimationfunctional}
\sup_{x\in\mathbb{R}}\left|\mathbb{P}\!\left\{\Phi\big(\gamma_{{Nk}}(\cdot)\big)
\leq x\right\} -\mathbb{P}\!\left\{\Phi\big(W_{{Nk}}(\cdot)\big)\leq x\right\}\right|
=O\left(N^{-1/4}(\log aN)^{3/4}\right).
\end{equation}
Assume further that the distribution of the r.v. $\Phi\big(V_{n}(F(\cdot))\big)$
 has a bounded density. Then, as $n\to\infty$,
\begin{equation}\label{estimationfunctional}
\sup_{x\in\mathbb{R}}\left|\mathbb{P}\!\left\{\Phi\big(\alpha_{n}(\cdot)\big)
\leq x\right\} -\mathbb{P}\!\left\{\Phi\big(V_n(\cdot)\big)\leq x\right\}\right|
=O\left(N^{-1/4}(\log aN)^{3/4}\right)\!.
\end{equation}
\end{corollary}
\noindent The proof of Corollary \ref{corol22} is postponed until Section \ref{proofc}.\\

\begin{remark}\label{remark33}
By the Borel-Cantelli Lemma and Theorem \ref{2eme theorem} we have
\begin{equation}
\sup_{0 \leq t \leq a
}|\gamma_{{Nk}}(t)-W_{{Nk}}(t)|{\stackrel{a.s.}{=}}
O\left(N^{-1/4}(\log aN)^{3/4}\right).
\end{equation}
\noindent Applying the Borel-Cantelli Lemma and Theorem \ref{3eme
theorem} we have
\begin{equation}
\sup_{0 \leq x \leq
Q_{k}(a)}|\alpha_n(x)-V_n(x)|{\stackrel{a.s.}{=}}
O\left(N^{-1/4}(\log aN)^{3/4}\right).
\end{equation}
\noindent For $a=1$, our results reduce to the results of
\cite{alybeirlanthorvart1984}.
\end{remark}

\begin{theorem} [\cite{alybeirlanthorvart1984}]
\label{theoreme1.1} Given the process $\{\alpha_{n}(x):0\leq x<\infty \}$
constructed from a sequence $U_{1},U_{2},\ldots $ of i.i.d random
variables of uniform low on $[0,1]$ and defined on a space of probability
eventually  enlarged version of $(\Omega ,A,\mathbb{P})$, there exists a sequence of 
Brownian bridges $B_{1},B_{2},\ldots ,$ defined on $(\Omega ,A,\mathbb{P})$
such that, for all $0\leq s<\infty $, if, 
\begin{equation}
V_n(x)=B_{N}(F_{k}(x))-\frac{1}{k}xf_{k}(x)\int_{0}^{\infty
}B_{N}(F_{k}(u))du,  \label{eq1.1 AlyBeirlHor84}
\end{equation}
then with probability $1$, 
\begin{equation}
\sup_{0\leq x<\infty }|\alpha_n(x)-\Gamma _{N}(x)|=O\left(
n^{-1/4}(\log n)^{3/4}\right) ,~~\mbox{as}~~n\rightarrow \infty .
\label{eq1.2 AlyBeirlHor84}
\end{equation}
\end{theorem}
\noindent Notice that the approximating Gaussian processes, for $k=1$, is given by 
\begin{equation*}
\Gamma _{n}(x)=B_{n}(x)+(1-x)\log(1-x)\int_{0}^{1}\frac{B_{n}(u)}{1-u}du.
\end{equation*}
\noindent A similar approximation is obtained for the
process of quantiles $\{\gamma _{N}(t):0\leq t\leq1, N\geq 1 \}$, in the special case where $N=n/k$. This
approximation is the following one.

\begin{theorem}[\cite{alybeirlanthorvart1984}]
\label{theoreme1.2} given the process $\{\gamma _{N}(t):0\leq t\leq 1, N\geq 1\}$
constructed from a sequence $U_{1},U_{2},\ldots $ of i.i.d random
variables of uniform low on $[0,1]$ and defined on a space of probability
eventually enlarged version of $(\Omega ,A,\mathbb{P})$, there exists a sequence of 
Brownian bridges $B_{1},B_{2},\ldots ,$ defined on $(\Omega ,A,\mathbb{P})$
such that, for all $0\leq t\leq 1$, if, 
\begin{equation}
W_{{Nk}}(t)=B_{N}(t)-\frac{1}{k}Q_{k}(t)f_{k}(Q_{k}(t))
\int_{0}^{1}B_{N}(u)dQ_{k}(u),  \label{eq2.1 AlyBeirlHor84}
\end{equation}
then with probability $1$, 
\begin{equation}
\sup_{0\leq t\leq 1}|\gamma _{N}(t)-W_{{Nk}}(t)|=O
\left( N^{-1/4}(\log N)^{3/4}\right), ~~\mbox{as}~~n\rightarrow \infty .
\label{eq2.2 AlyBeirlHor84}
\end{equation}
\end{theorem}

\section{Functional limit laws for the increments of the quantile process of
non-overlapping $k$-spacings processes}\label{sectionf}

\noindent Hereafter, our study of the process of quantile
of $k$-spacings $\{\gamma _{N}(t):0\leq t\leq 1\}$,  will be restricted  to the
particular case where $N=n/k$. The increments of the process $\{\gamma
_{N}(t):0\leq t\leq 1\}$ are defined in the following way: for all $0<h<1$
and $n\geq k$ 
\begin{equation}
\eta _{N}(h,t;s)=\gamma _{N}(t+sh)-\gamma _{N}(t),~\mbox{for}~0\leq t\leq
1~\mbox{and}~s\in \mathbb{R}.  \label{eq0}
\end{equation}
\noindent In this section, we are interested in the study of the fluctuations of the process of quantile. We are
going to give a functional law of the iterated logarithm for the
increments $\{\eta _{N}(h,t;s): 0\leq t\leq
1, -\infty<s <\infty\}$. \cite{Dindar1997} established  the functional law of the iterated
logarithm for the increments of the processes associated to the uniform
spacing. Our aim, is to provide analogous results for the $k$-spacing quantile process. In particular, we are going to prove that the results obtained by \cite{Dindar1997}  for
the process of quantile of the uniform spacings, are still valid in the
general case of the uniform non-overlapping $k$-spacings by preserving
the same hypotheses.

\vskip5pt \noindent In the sequel, we fix the
following notation: Let $B[0,1]$ denote the set of bounded functions
defined on $[0,1].$ $B[0,1]$ is endowed with the topology induced by the
norm sup $$\Vert f\Vert =\sup_{0\leq s\leq 1}|f(s)|.$$
Let $AC[0,1]$ be the set of all absolutely continuous function on $[0,1]$,
and $\dot{f}(s)=df(s)/ds$ the Lebesgue derivative of $f(\cdot)\in AC[0,1]$. For any 
$c\geq 0$, we denote by 
\begin{equation*}
\mathcal{S}_{c}=\left\{ f(\cdot)\in AC[0,1]:f(0)=0~\mbox{et}~\int_{0}^{1}\dot{f}
^{2}(s)ds\leq c\right\}
\end{equation*}
the Strassen set (see e.g., \cite{Strassen1964}). For any subset $A\subset
AC[0,1]$ and $\varepsilon >0$, let $$A^{\varepsilon }=\{f(\cdot)\in B[0,1]:\exists
g(\cdot)\in A,\Vert f-g\Vert <\varepsilon \}.$$ The Hausdorff distance between  the sets $A$
and $B$ is defined to be 
\begin{equation*}
\boldsymbol{\Delta} (A,B)=\inf \{\varepsilon >0:A\subseteq B^{\varepsilon }~\mbox{and}
~B\subseteq A^{\varepsilon }\}.
\end{equation*}
Consider a sequence of constants $\{h_{N}:N\geq 1\}$ satisfying the
following assumptions

\begin{description}
\item[(H.1)] $h_{N}\downarrow 0~\mbox{and}~Nh_{N}\uparrow \infty~ 
\mbox{where}~ N\uparrow \infty~\mbox{and}~ 0<h_{N}<1$,

\item[(H.2)] $\log(1/h_{N})/\log_2 N \rightarrow c\in [0,\infty], ~
\mbox{where}~N\rightarrow \infty$,

\item[(H.3)] $(1/h_{N})=o\big(N^{1/2}(\log_2 N)/(\log N)^{3/2}\big),~
\mbox{where}~N\rightarrow \infty$.
\end{description}
 Let $\overset{\mathbb{P}}{=}$ denotes equality in probability.
\noindent Let $$\mathcal{\beta}_{N}=\big(2h_{N}\big\{\log(1/h_{N})+\log_2 N
\big\}\big)^{-1/2}, ~~\mbox{for}~~N\geq 3.$$ 
The following theorem constitute our main results of this section.
\begin{theorem}\label{theoreme increment m spacings}
Assume that $\{h_{N} : n \geq  1\}$ fulfills the limiting conditions,\textbf{(H.1)}, \textbf{(H.2)}, \textbf{(H.3)}. Then we have 
\begin{equation}
\lim_{N\rightarrow \infty}\mathbf{\Delta}\Big(\big\{\mathcal{\beta}_{N}\eta_{N}(h_{N},t;s)~:~0\leq
t\leq1-h_{N}\big\},\mathcal{S}_{\frac{c}{c+1}}\Big)\overset{\mathbb{P}}{=}0.
\end{equation}
\end{theorem}
\noindent The proof of Theorem \ref{theoreme increment m spacings} is postponed until Section \ref{proofc}.
\subsection{Limit laws of functional of the $k$-spacings process}

\vskip5pt \noindent Let $\Phi :B_{0}[0,1]\rightarrow \mathbb{R}$ be a
functional defined on a subset   $B_{0}[0,1]$ of $B[0,1]$ such that

\begin{description}
\item[(i)] $\beta_{N} \eta_{N}(h_{N},t;\cdot)\in B_0[0,1]~~ \forall~~ 0\leq
t\leq 1-h_{N}$,

\item[(ii)] $\mathcal{S}_c\subseteq B_0[0,1]~~\forall~~c>0$,

\item[(iii)] $\Phi $ is continuous for the norm sup in $B_{0}[0,1]$.
\end{description}

\noindent For $f(\cdot)\in B_{0}[0,1]$, examples of such functional are given by the following 

\begin{itemize}
\item $\Phi_1(f)=|f(1)|$,

\item $\Phi_2(f)=\sup_{0\leq s\leq 1}|f(s)|$,

\item $\Phi_3(f)=\pm\int_0^1f(u)dv(u)$,
where we suppose that $v(\cdot)$ has bounded variation with $v(1)=0$.
\end{itemize}

\begin{remark}
We notice that $\Phi _{3}(f)$ is not defined for all $f(\cdot)\in B[0,1]$, which
justify the introduction of $B_{0}[0,1]$.
\end{remark}

\vskip5pt Let us consider the process defined by $\{\Phi (\beta _{N}\eta
_{N}(h_{N},t;\cdot)):0\leq t\leq 1-h_{N}\}$. We have the following results.

\begin{corollary}
\label{corollaire1} Assume that $\{h_{N} : N \geq  1\}$ fulfills the limiting conditions,\textbf{(H.1)}, \textbf{(H.2)}, \textbf{(H.3)}. Then we have 
\begin{equation}
\lim_{n\rightarrow \infty }\sup_{0\leq t\leq 1-h_{N}}\Phi (\beta _{N}\eta
_{N}(h_{N},t;\cdot))\overset{\mathbb{P}}{=}\sup_{f\in \mathcal{S}_{\frac{c}{c+1}}}\Phi
(f).
\end{equation}
\end{corollary}
\noindent The proof of  Corollary \ref{corollaire1} is postponed until Section \ref{proofc}.
\\ By applying Corollary \ref{corollaire1} to the 
functionals $\Phi $  given above in example, we obtain the following corollary.

\begin{corollary}
\label{corollaire2} Assume that $\{h_{N} : N \geq  1\}$ fulfills the limiting conditions, \textbf{(H.1)}, \textbf{(H.2)}, \textbf{(H.3)}. Then we have \begin{eqnarray}
\lim_{n\rightarrow \infty }\left[ \sup_{0\leq t\leq 1-h_{N}}\beta
_{N}|\gamma _{N}(t+h_{N})-\gamma _{N}(t)|\right] &\overset{\mathbb{P}}{=}&\left( \frac{c
}{c+1}\right) ^{1/2},  \label{eq1 corrolaire2}
\\
\lim_{n\rightarrow \infty }\left[ \sup_{0\leq t\leq 1-h_{N}}\sup_{0\leq
s\leq 1}\beta _{N}|\gamma _{N}(t+h_{N}s)-\gamma _{N}(t)|\right] &\overset{\mathbb{P}}{=
}&\left( \frac{c}{c+1}\right) ^{1/2},  \label{eq2 corrolaire2}
\\
\lim_{n\rightarrow \infty }\left[ \sup_{0\leq t\leq 1-h_{N}}\pm
\int_{0}^{1}\beta _{N}\eta _{N}(h_{N},t;s)dv(s)\right] &\overset{\mathbb{P}}{=}&\left( 
\frac{c}{c+1}\int_{0}^{1}v^{2}(s)ds\right) ^{1/2}.  \label{eq3 corrolaire2}
\end{eqnarray}
\end{corollary}
\noindent The proof of  Corollary \ref{corollaire2} is postponed until Section \ref{proofc}.

\begin{remark}
The results given above for the increments of the process of quantile of the
disjoint $k$-spacings, and the results achieved in \cite{Dindar1997} for the
empirical process of the uniform spacings, remain valid for the increments
of the empirical process $\{\alpha_{n}(x ):0\leq x< \infty, n\geq 1\}$.
\end{remark}

\section{Proof}\label{proofc}
This section is devoted to the proofs of our results. The previously defined notation continues to be used in the following.

\subsection*{Proof of Theorem \ref{1er theorem}.}\label{th1} Consider the
sequence $\xi_i=F_{k}(Y_i), i=1,2,\ldots,$ of i.i.d. $U[0,1]$
$r.v^{,}s$ and construct the corresponding uniform quantile
process defined by
\begin{equation}\label{U_N(t)}
U_N(t) = N^{1/2}(t- F_{k}(K_N(t))),
\end{equation}
where $Y_i$ and $K_N(\cdot)$ are defined by (\ref{Y_i}) and
(\ref{K_N(t)}) successively. A simple application of Theorem
$1.1$ of \cite{CsCsHM1986} with $a=d/n$ and
$x=\varepsilon\lambda^{-1}\log aN$, we can find a sequence of
Brownian bridges $\{B_N^{(2)}(t) : 0\leq t \leq\ 1\}$, such that for
all $\varepsilon
>0$ we have
\begin{equation}\label{Masonvanzwet}
\mathbb{P}\left( \sup_{0 \leq t \leq a}|U_N(t)-B_N^{(2)}(t)| \geq
A_{4}N^{-1/2}(\log aN) \right) \leq B_{4}N^{-\varepsilon},
\end{equation}
 where $A_{4}, B_{4}$ are positive constants depending on $\varepsilon$ and $a$.
 Furthermore, we have for all $0 \leq a \leq 1$,
\begin{equation}
 \mathbb{P}\left( \sup_{0 \leq t \leq a }|B_N^{(2)}(t)| > x\right) \leq 2e^{-2x^2},~~ x \geq 0.\\
\end{equation}
\noindent The last inequality together with (\ref{Masonvanzwet})
implies that
\begin{equation}\label{P(sup_(0,a)|U_N(t)|)}
\mathbb{P}\left( \sup_{0\leq t\leq
a}|U_N(t)|\geq\left(\frac{1}{2}\varepsilon (\log aN)
\right)^{1/2}+A_{4}N^{-1/2}(\log aN)\right)\leq
\left(2+B_{4}\right)N^{-\varepsilon}.
\end{equation}
\noindent We will prove in the next lemma that $\{U_N(t): 0\leq t\leq 1, N\geq 1\}$, as
defined in (\ref{U_N(t)}), can be approximated by $\{B^{(1)}_N(t): 0\leq t\leq 1, N\geq 1\}$ as
well.
\vskip5pt
\begin{lemma}\label{1er lemme}
There exists a sequence of Gaussian processes $\{B^{(1)}_N(t): 0\leq t\leq 1, N\geq 1\}$, such that, for all $\varepsilon>0$, we have
\begin{equation}
\mathbb{P}\left( \sup_{0 \leq t \leq a}|U_N(t) - B_N^{(1)}(t)| \geq
A_{5}N^{-1/2}(\log aN)^{3/4}\right)\leq B_{5}N^{-\varepsilon},
\end{equation}
where $A_{5}$ and $B_{5}$  are positive constants.
\end{lemma}
\vskip5pt
\noindent {\bf{Proof of Lemma \ref{1er lemme}}.} Let
$\xi_{1,N},\ldots,\xi_{N,N}$  denote the order statistics of
$\xi_{1},\ldots,\xi_{N}$. \noindent By Theorem \ref{Masonvanzwet
theorem} and the fact that
$$\beta_N(Q_{k}(\xi_{i,N}))=U_N\left(\frac{i}{N}\right),$$ we have, for each
$0<a\leq 1$,
\begin{equation}
\mathbb{P}\left (\max_{0 \leq i \leq aN}\left|U_N\left(\frac{i}{N}\right)-
B^{(1)}_N\left(\xi_{i,N}\right)\right|>
\mathcal{A}N^{-1/2}\left(\log aN\right) \right ) \leq
\mathcal{B}N^{-\varepsilon}.
\end{equation}
\noindent On the other hand, from (\ref{P(sup_(0,a)|U_N(t)|)}) we
have
\begin{equation}\label{equ3.5}
\mathbb{P}\left (\max_{0 \leq i \leq aN}\left|
\frac{i}{N}-\xi_{i,N}\right|\geq
 N^{-1/2}\left(\frac{\varepsilon}{2}(\log aN)\right)^{1/2}+
 A_{4}N^{-1}\left(\log aN\right)\right)\leq
(2+B_{4})N^{-\varepsilon}.
\end{equation}
\noindent An application of Lemma  $1.2.1$ in connection with Lemma $1.4.1$ of
\cite{csorgorevesz1981} allow us to write
\begin{eqnarray*}
 \mathbb{P}\left( \sup_{0 \leq i \leq N-N^{1/2}(\log aN)}\sup_{0 \leq s \leq
 N^{-1/2}(\log aN)}\left|B^{(1)}_N\left(\frac{i}{N}+s\right)-B^{(1)}_N\left(\frac{i}{N}\right)\right|
\right.\\\left.> A_{6}N^{-1/4}(\log aN)^{3/4} \right )\leq
B_{6}N^{-\varepsilon},
 \end{eqnarray*}
\noindent This, when combined with (\ref{equ3.5}), implies that
 \begin{equation}
 \mathbb{P}\left(\max_{0 \leq i \leq
  aN}|B^{(1)}_N\left(\frac{i}{N}\right)-B^{(1)}_N\left(\xi_{i,N}\right)|> A_{7}N^{-1/4}(\log
  aN)^{3/4}\right )\leq B_{7}N^{-\varepsilon}.
  \end{equation}
\noindent Lemma \ref{1er lemme} follows from the fact that
\begin{equation}
\left|U_N(t)-U_N\left(\frac{i}{N}\right)\right|\leq
N^{-1/2}~~\mbox{for}~~ \frac{i-1}{N}< t <\frac{i}{N}.
\end{equation}

\hfill$\blacksquare$\vskip5pt\noindent We return now to the proof
of Theorem \ref{1er theorem}. Following
\cite{alybeirlanthorvart1984}, we have
\begin{eqnarray}
\sup_{0<t<\infty}F_{k}(t)(1-F_{k}(t))
\frac{|(f_{k})^{'}(t)|}{(f_{k})^2(t)} &\leq& \gamma,
\end{eqnarray}
together with
\begin{equation}
\lim_{t\rightarrow \infty}F_{k}(t)(1-F_{k}(t))
\frac{|(f_{k})^{'}(t)|}{(f_{k})^2(t)}=1,
\end{equation}
\begin{equation}
\lim_{t\rightarrow 0}F_{k}(t)(1-F_{k}(t))
\frac{|(f_{k})^{'}(t)|}{(f_{k})^2(t)}=1,
\end{equation}
for some $$\gamma=\gamma(k) < \infty.$$ By the
mean value theorem, we readily obtain
\begin{equation}\label{K_N-U_N}
\kappa_N(t)-U_N(t)=U_N(t)\left(\frac{f_{k}(Q_{k}(t))}{f_{k}(Q_{k}(\theta_{t,N}))}-1\right),
\end{equation}
\noindent for some $\theta_{t,N}$ such that
$$|\theta_{t,N}-t|<N^{-1/2}|U_N(t)|.$$ In Theorem $1.5.1$ in
\cite{csorgo1983}, it is proved that
\begin{eqnarray}\label{csorgo93}
\lefteqn{ \mathbb{P}\left(\sup_{c\leq t \leq 1-c}\left|
\frac{f_{k}(Q_{k}(t))}{f_{k}(Q_{k}(\theta_{t,N}))}-1\right|>\delta
\right)}\nonumber\\&\leq&
4([\gamma]+1)\{\exp(-N{\rm ch}((1+\delta)^{1/2([\gamma]+1)}))
\nonumber\\&&+\exp(-N{\rm ch}((1+\delta)^{-1/2([\gamma]+1)}))\},
\end{eqnarray}
for all $\delta >0, 0<c<1$ and $N\geq 1$, where $$h(x)=x+\log
(1/x)-1, ~\mbox{for}~x~x>0.$$
Moreover, there exist a $\delta_0> 0$ such that
\begin{equation}
h((1+\delta)^{\mp1/2([\gamma]+1)})\geq
\frac{1}{8}([\gamma]+1)^2\delta^2,~\mbox{for}~0<\delta <\delta_0.
\end{equation}
Let $$\delta_N:=(8\varepsilon)^{1/2}\left([\gamma]+1)^{-1}N^{-1/4}(\log
aN\right)^{1/2},$$ and  $$C^{(1)}:=C_N^{(1)}:=N^{-1/2}.$$  By the above
inequality and (\ref{csorgo93}) we obtain that, for $N$
sufficiently large, that
\begin{equation}\label{sup_(c<t<1-c)}
 \mathbb{P}\left(\sup_{C_N^{(1)}\leq t \leq 1-{C_N^{(1)}}}\left|
\frac{f_{k}(Q_{k}(t))}{f_{k}(Q_{k}(\theta_{t,N}))}-1\right|>\delta_N
\right)\leq 8([\gamma]+1)N^{-\varepsilon}.
\end{equation}

\noindent Combining (\ref{K_N-U_N}), (\ref{P(sup_(0,a)|U_N(t)|)})
and (\ref{sup_(c<t<1-c)}), we obtain that, for $N$ sufficiently
large, 
\begin{equation}\label{sup_(0<t<a-c)(K_N-U_N)}
 \mathbb{P}\left(\sup_{C_N^{(1)}\leq t \leq a-C_N^{(1)}}\left|\kappa_N(t)- U_N(t) \right|>
 A_{8}N^{-1/4}(\log aN)^{3/4}\right)\leq B_{8}N^{-\varepsilon}.
\end{equation}
To complete the proof of Theorem \ref{1er theorem}, we replace
$\log N$ in the proof of the Theorem B of
\cite{alybeirlanthorvart1984} by $(\log aN)$. We omit the details, which essentially repeats, more or the less verbatim,
the same arguments.
\hfill$\blacksquare$\\
The following technical  Lemmas \ref{2eme lemme} and \ref{lemma 3} will be instrumental in the proof of Theorem \ref{2eme theorem}.
\vskip5pt
\begin{lemma}\label{2eme lemme}
We have, for each $\varepsilon>0$, and all $n\geq m$ sufficiently
large
\begin{equation}
\mathbb{P}\left(\left|N^{1/2}\left(\frac{T_N}{{Nk}}-1\right)-
\frac{1}{k}\int_0^{\infty}tdB^{(1)}_N\left(F_{k}(t)\right)\right|
>A_{9}N^{-1/2}\left(\log aN\right)^{2}\right)\leq
B_{9}N^{-\varepsilon}.
\end{equation}
 where $A_{9}=A_9(\varepsilon)=4(1/2+\varepsilon)
 \mathcal{A}$ and $B_{9}=8\sqrt{2}+\mathcal{B}$
denote positive constants.
\end{lemma}
\vskip5pt
\noindent {\bf{Proof of Lemma \ref{2eme lemme}.}}\vskip5pt
\noindent It is readily checked that,
\begin{equation}
\frac{T_N}{{Nk}}=\frac{1}{{Nk}}\sum_{i=1}^{N}Y_i=\frac{1}{k}\int_0^{\infty}tdG_N(t)~~
\mbox{and}~~\int_0^{\infty}tdF_{k}(t)=k.
\end{equation}
\noindent From which we obtain readily that
\begin{equation}
N^{1/2}\left(\frac{T_N}{N}-k\right)=\int_0^{\infty}td\beta_N(t)=-\int_0^{\infty}\beta_N(t)dt.
\end{equation}
\noindent Let $\lambda_N$ be a sequence of positive numbers. Making use of the triangle inequality, we infer that 
\begin{eqnarray*}
\left|\int_0^{\infty}\beta_N(t)dt-\int_0^{\infty}B^{(1)}_N(F_{k}(t))dt\right|
&\leq&\int_0^{\lambda_N}\left|\beta_N(t)-B^{(1)}_N(F_{k}(t))\right|dt
\\&&+\int_{\lambda_N}^{\infty}|B^{(1)}_N(F_{k}(t))|dt\\&&+ \int_{\lambda_N}^{\infty}|\beta_N(t)|dt.
\end{eqnarray*}
\noindent Let us recall the following well known properties 
\begin{eqnarray}
\mathbb{E}\left(\beta_N(t)\right)&=&\mathbb{E}\left(B^{(1)}_N(F_{k}(t))\right)=0,
\\
{\rm Var}\left(\beta_N(t)\right)&=&\mathbb{E}\left[\left(\beta_N(t)\right)^2\right]=F_{k}(t)(1-F_{k}(t)),
\end{eqnarray}
and
\begin{equation}
{\rm Var}\left(B^{(1)}_N(F_{k}(t)\right)=\mathbb{E}\left[\left(B^{(1)}_N(F_{k}(t)\right)^2\right]
=F_{k}(t)(1-F_{k}(t)).
\end{equation}
\noindent  Making use of the Fubini theorem's, in combination with Cauchy-Schwartz inequality implies that 

\begin{eqnarray}\label{cauchy-sh1}
\mathbb{E}\int_{\lambda_N}^{\infty}|\beta_N(t)|dt&=&\int_{\lambda_N}^{\infty}\mathbb{E}|\beta_N(t)|dt\nonumber\\&\leq&
\int_{\lambda_N}^{\infty}(F_{k}(t)(1-F_{k}(t)))^{1/2}dt,
\end{eqnarray}

\noindent and similarly 
\begin{eqnarray}\label{cauchy-sh2}
\mathbb{E}\int_{\lambda_N}^{\infty}|B^{(1)}_N(F_{k}(t))|dt&=&\int_{\lambda_N}^{\infty}\mathbb{E}|B^{(1)}_N(F_{k}(t))|dt
\nonumber\\&\leq&
\int_{\lambda_N}^{\infty}(F_{k}(t)(1-F_{k}(t)))^{1/2}dt.
\end{eqnarray}
\noindent By \cite{alybeirlanthorvart1984}, there exists $t_0>0$
such that
\begin{equation}\label{aaa1}
1-F_{k}(t)\leq2\exp\left(-\frac{t}{2}\right), ~~\mbox{if}~~
t\geq t_0.
\end{equation}
\noindent Hence, provided that $\lambda_N\geq t_0$, by
(\ref{aaa1}) and the fact that
\begin{equation}
F_{k}(t)\leq 1 ~~\mbox{for all}~~ t>0,
\end{equation}
\noindent the left hand sides of (\ref{cauchy-sh1}) and
(\ref{cauchy-sh2}) are bounded above by
$4\sqrt{2}\exp(-\lambda_N/4)$. Indeed, we have 
\begin{eqnarray*}
\mathbb{E}\left(|B^{(1)}_N(F_{k}(t))|\right)&\leq&
(F_{k}(t)(1-F_{k}(t)))^{1/2}\\&\leq&\sqrt{2}\exp(-t/4),
\end{eqnarray*}
and by using (\ref{cauchy-sh2}), we infer that
\begin{eqnarray*}
\mathbb{E}\left(\int_{\lambda_N}^{\infty}|B^{(1)}_N(F_{k}(t))|dt\right)
&\leq&\sqrt{2}\int_{\lambda_N}^{\infty}\exp(-t/4)dt\\&=&
4\sqrt{2}\exp(-\lambda_N/4).
\end{eqnarray*}
By using similar arguments shows likewise that
\begin{equation*}
\mathbb{E}\left(\int_{\lambda_N}^{\infty}|\beta_N(t)|dt\right)\leq
4\sqrt{2}\exp(-\lambda_N/4).
\end{equation*}
\noindent By choosing $\lambda_N=4(\frac{1}{2}+\varepsilon)(\log
aN)$, and applying Markov inequality, we infer that
\begin{equation}
\mathbb{P}\left(\int_{4(\frac{1}{2}+\varepsilon)(\log
aN)}^{\infty}|\beta_N(t)|dt>
a^{-(1/2+\varepsilon)}N^{-1/2}\right)\leq
4\sqrt{2}N^{-\varepsilon},
\end{equation}
\noindent and
\begin{equation}
\mathbb{P}\left(\int_{4(\frac{1}{2}+\varepsilon)(\log
aN)}^{\infty}|B^{(1)}_N(F_{k}(t))|dt>
a^{-(1/2+\varepsilon)}N^{-1/2}\right)\leq
4\sqrt{2}N^{-\varepsilon}.
\end{equation}
\noindent   An application of Theorem \ref{Masonvanzwet theorem} shows that, 
\begin{equation}\label{eq333}
\mathbb{P}\left(\int_0^{\lambda_N}\left|\beta_N(t)-B^{(1)}_N(F_{k}(t))\right|dt>\lambda_N\mathcal{A}N^{-1/2}(\log
aN)\right)\leq \mathcal{B}N^{-\varepsilon}.
\end{equation}
\noindent More precisely, we have the inequality
\begin{eqnarray*}
\int_0^{\lambda_N}\left|\beta_N(t)-B^{(1)}_N(F_{k}(t))\right|dt&\leq&
\sup_{0\leq x\leq
Q_{k}(a)}\left|\beta_N(t)-B^{(1)}_N(F_{k}(t))\right|\int_0^{\lambda_N}dt
\\&=&\lambda_N\sup_{0\leq x\leq
Q_{k}(a)}\left|\beta_N(t)-B^{(1)}_N(F_{k}(t))\right|.
\end{eqnarray*}
Making use of  Theorem \ref{Masonvanzwet theorem}, we obtain 
\begin{eqnarray*}
\lefteqn{\mathbb{P}\left(\int_0^{\lambda_N}\left|\beta_N(t)-B^{(1)}_N(F_{k}(t))\right|dt>\lambda_N\mathcal{A}N^{-1/2}(\log
aN)\right)}\\&\leq&\mathbb{P}\left(\lambda_N\sup_{0\leq x\leq
Q_{k}(a)}\left|\beta_N(t)-B^{(1)}_N(F_{k}(t))\right|>
\lambda_N\mathcal{A}N^{-1/2}(\log aN)\right)\\&=&\mathbb{P}\left(
\sup_{0\leq x\leq
Q_{k}(a)}\left|\beta_N(t)-B^{(1)}_N(F_{k}(t))\right|>\mathcal{A}N^{-1/2}(\log
aN)\right)\\&\leq& \mathcal{B}N^{-\varepsilon}.
\end{eqnarray*}
Let $$\Lambda_1=2a^{-(1/2+\varepsilon)}N^{-1/2},$$ and
$$\Lambda_2=\lambda_N\mathcal{A}N^{-1/2}(\log
aN)=4(1/2+\varepsilon)\mathcal{A}N^{-1/2}(\log aN)^{2}.$$ From this we infer that we have  
\begin{equation*}
\Lambda_1+\Lambda_2=4(1/2+\varepsilon)\mathcal{A}N^{-1/2}(\log
aN)^{2}\left(1+o(1)\right).
\end{equation*}

\noindent Lemma \ref{2eme lemme} now follows by combining the
above three inequalities (\ref{cauchy-sh1}), (\ref{cauchy-sh2}) and (\ref{eq333}). We have 
\begin{eqnarray*}
\lefteqn{\mathbb{P}\left(\left|\int_0^\infty\left(\beta_N(t)-B^{(1)}_N(F_{k}(t))\right)dt\right|>
\Lambda_1+\Lambda_2\right)}\\&\leq&\mathbb{P}\left(\left|\int_0^\lambda\left(\beta_N(t)-B^{(1)}_N(F_{k}(t))\right)dt\right|>
4(1/2+\varepsilon)\mathcal{A}N^{-1/2}(\log_+
aN)^{2}\right)\\&&+\mathbb{P}\left(\left|\int_{\lambda}^\infty\left(B^{(1)}_N(F_{k}(t))\right)dt\right|>
a^{-(1/2+\varepsilon)}N^{-1/2}\right)\\&&+\mathbb{P}\left(\left|\int_{\lambda}^\infty\left(\beta_N(t)\right)dt\right|>
a^{-(1/2+\varepsilon)}N^{-1/2}\right)\\&\leq&
4\sqrt{2}N^{-\varepsilon}+ 4\sqrt{2}N^{-\varepsilon}+\mathcal{B}.
\end{eqnarray*}
By choosing $A_9=A_9(\varepsilon)=4(1/2+\varepsilon)\mathcal{A}$
and $B_9=8\sqrt{2}+\mathcal{B}$,  the proof of Lemma \ref{2eme
lemme} is  completed. \hfill$\blacksquare$\vskip5pt
\begin{lemma}\label{lemma 3}
For each $\varepsilon >0$ and $n\geq k$, we have, uniformly over
$0 \leq a \leq 1$
\begin{eqnarray*} \mathbb{P}\left(\sup_{0 \leq x\leq
Q_{k}(a)}\left|B^{(1)}_N\left(F_{k}\left(x\frac{T_N}{{Nk}}\right)\right)
-B^{(1)}_N(F_{k}(x))\right|\right.\\\left.>A_{10}N^{-1/4}(\log
aN)^{3/4}\right)&\leq& B_{10}N^{-\varepsilon},
\end{eqnarray*}
\noindent where $A_{10}$ and  $ B_{10}$ are positive constants.
\end{lemma}
\vskip5pt
\noindent {\bf{Proof of lemma \ref{lemma 3}.}} The random variable
$\int_{0}^{\infty}B^{(1)}_N(F_{k}(t))dt$ has a normal
distribution, with expectation $0$ and finite variance, given by
\begin{equation}
\sigma^2_1=\mathbb{E}\left\{\left(\int_{0}^{\infty}B^{(1)}_N(F_{k}(t))dt\right)^2\right\}<
\infty.
\end{equation}
\noindent Hence, we have the inequality
\begin{equation}\label{inequality (3.57)}
\mathbb{P}\left(\frac{1}{\sigma_1}\left|\int_{0}^{\infty}B^{(1)}_N(F_{k}(t))dt\right|>(2\varepsilon\log
aN)^{1/2}\right)\leq 2N^{-\varepsilon}.
\end{equation}
 \noindent This inequality and Lemma \ref{2eme lemme} imply that
\begin{equation}\label{inequality (3.6)}
\mathbb{P}\left(\left|\frac{T_N}{{Nk}}-1\right|>A_{11}N^{-1/2}(\log
aN)^{1/2}\right)\leq B_{11}N^{-\varepsilon},
\end{equation}
where
$$A_{11}=A_{11}(\varepsilon)=\left(2k^{-2}\sigma_1^2\varepsilon
\right)^{1/2} ~~\mbox{and}~~ B_{11}=2+B_{9}.$$ More precisely, we have 
$$A_{9}N^{-1/2}(\log aN)^{1/2}+(2k^{-2}\sigma_1^2\varepsilon \log
aN)^{1/2}=(2k^{-2}\sigma_1^2\varepsilon)^{1/2}( \log
aN)^{1/2}(1+o(1)).$$ So the probability (\ref{inequality (3.6)})
is the same as
$$\mathbb{P}\left(\left|N^{1/2}\left(\frac{T_N}{{Nk}}-1\right)\right|>A_{9}N^{-1/2}(\log
aN)^{1/2}+(2k^{-2}\sigma_1^2\varepsilon \log aN)^{1/2}\right).$$
\noindent By combination of Lemma \ref{2eme lemme} and inequality
(\ref{inequality (3.57)}), it follows that 
\begin{eqnarray*}
\lefteqn{\mathbb{P}\left(\left|N^{1/2}\left(\frac{T_N}{{Nk}}-1\right)\right|>A_{9}N^{-1/2}(\log
aN)^{1/2}+(2k^{-2}\sigma_1^2\varepsilon \log
aN)^{1/2}\right)}\\&\leq&\mathbb{P}\left(\left|N^{1/2}\left(\frac{T_N}{{Nk}}-1\right)-
\frac{1}{k}\int_{0}^{\infty}B^{(1)}_N(F_{k}(t))dt\right|>A_{9}N^{-1/2}(\log
aN)^{1/2}\right)\\&&
+\mathbb{P}\left(\left|\frac{1}{k}\int_{0}^{\infty}B^{(1)}_N(F_{k}(t))dt\right|>(2k^{-2}\sigma_1^2\varepsilon
\log
aN)^{1/2}\right)\\&=&\mathbb{P}\left(\left|N^{1/2}\left(\frac{T_N}{{Nk}}-1\right)-
\frac{1}{k}\int_{0}^{\infty}B^{(1)}_N(F_{k}(t))dt\right|>A_{9}N^{-1/2}(\log
aN)^{1/2}\right)\\&&
+\mathbb{P}\left(\left|\frac{1}{\sigma_1}\int_{0}^{\infty}B^{(1)}_N(F_{k}(t))dt\right|>(2\varepsilon
\log aN)^{1/2}\right)\\&\leq& (B_9+2)N^{-\varepsilon}.
\end{eqnarray*}
\vskip5pt \noindent 
 By  Taylor expansions, we readily obtain 
\begin{equation}\label{inequality (3.7)}
\left|F_{k}\left(x\frac{T_N}{{Nk}}\right)-F_{k}(x)\right|=xf_{k}(x_N)\left|\frac{T_N}{{Nk}}-1\right|,
\end{equation}
\noindent where $$|x_N-x|\leq x\left|\frac{T_N}{{Nk}}-1\right|.$$
\noindent Let $0<\delta<1$ and define $A_N(\delta)$ by
\begin{equation}
A_N(\delta)=\left\{\omega: \left|\frac{T_N}{{Nk}}-1\right|\leq
\delta\right\}.
\end{equation}
\noindent Now, by choosing $N$ sufficiently large so that
$A_{11}N^{-1/2}(\log aN)^{1/2}\leq \delta$, and using
(\ref{inequality (3.6)}) we get that
$$\mathbb{P}\left(A^c_N(\delta)\right)\leq B_{11}N^{-\varepsilon}.$$ In
addition, we have for each $x_N \in A_N(\delta)$,
\begin{equation}
xf_{k}(x_N)\leq
\frac{(1+\delta)^{m-1}}{\Gamma(k)}x^ke^{-(1-\delta)x},
\end{equation}
which is bounded on $[0,\infty)$. Now, we let 
\begin{equation}
A_{12}=A_{11}.\sup_{0 \leq x\leq
Q_{k}(a)}\frac{(1+\delta)^{k-1}}{\Gamma(k)}x^ke^{-(1-\delta)x}.
\end{equation}
\noindent Recall the following elementary fact  $$\mathbb{P}(A)\leq \mathbb{P}(B^c)+\mathbb{P}(A\cap B),$$then, for large enough $N$, we infer that we have 
\begin{eqnarray}\label{inequality (3.8)}
\lefteqn{\mathbb{P}\left(\sup_{0 \leq x\leq
Q_{k}(a)}\left|F_{k}\left(x\frac{T_N}{{Nk}}\right)-F_{k}(x)\right|>A_{12}N^{-1/2}(\log
aN)^{1/2}\right)}\nonumber\\&\leq&
\mathbb{P}\left(A^c_N(\delta)\right)\nonumber\\&&+
\mathbb{P}\left(A_N(\delta)~~\mbox{and}\left\{\sup_{0 \leq x\leq
Q_{k}(a)}\left|F_{k}\left(x\frac{T_N}{{Nk}}\right)-F_{k}(x)\right|>A_{12}N^{-1/2}(\log
aN)^{1/2}\right\}\right)\nonumber\\&\leq&
\mathbb{P}\left(A^c_N(\delta)\right)\nonumber\\&&+\mathbb{P}\left(A_N(\delta)~~\mbox{and}
\left\{\sup_{0 \leq x\leq
Q_{k}(a)}xF_{k}(x_N)\left|\frac{T_N}{{Nk}}-1\right|>A_{12}N^{-1/2}(\log
aN)^{1/2}\right\}\right) \nonumber\\&\leq&
B_{11}N^{-\varepsilon}+\mathbb{P}\left(A_N(\delta)~~\mbox{and}
\left\{\left|\frac{T_N}{{Nk}}-1\right|>A_{11}N^{-1/2}(\log
aN)^{1/2}\right\}\right)\nonumber\\&\leq&
B_{11}N^{-\varepsilon}.
\end{eqnarray}
\noindent Now, (\ref{inequality (3.8)}) when combined with Lemma
$1.1.1$ of \cite{csorgorevesz1981} implies that
\begin{eqnarray}
\lefteqn{\mathbb{P}\left(\sup_{0 \leq x\leq
Q_{k}(a)}\left|B^{(1)}_N\left(F_{k}\left(x\frac{T_N}{{Nk}}\right)\right)
-B^{(1)}_N(F_{k}(x))\right|>A_{10}N^{-1/4}(\log
aN)^{3/4}\right)}\nonumber\\&=&\mathbb{P}\left(\sup_{0 \leq x\leq
Q_{k}(a)}\left|B^{(1)}_N\left(F_{k}(x)+F_{k}
\left(x\frac{T_N}{{Nk}}\right)-F_{k}(x)\right)-B^{(1)}_N(F_{k}(x))\right|
\right. \nonumber \\& & \left. > A_{10}N^{-1/4}(\log
aN)^{3/4}\right)\nonumber\\ &\leq& \mathbb{P}\left(\sup_{0\leq t\leq
1-A_{12}N^{-1/2}(\log aN)^{1/2}}\sup_{0\leq s\leq
A_{12}N^{-1/2}(\log aN)^{1/2}}\left|B^{(1)}_N\left(t+s\right)
-B^{(1)}_N(t)\right|\right. \nonumber \\& & \left. >
\frac{A_{10}}{\sqrt{A_{12}}}(\log
aN)^{1/2}\left(A_{12}N^{-1/2}(\log
aN)^{1/2}\right)^{1/2}\right)+B_{11}N^{-\varepsilon}\nonumber\\&\leq&
B_{10}N^{-\varepsilon}.
\end{eqnarray}
\noindent This completes the proof of Lemma \ref{lemma
3}.\hfill$\blacksquare$
\subsection*{Proof of Theorem \ref{2eme
theorem}.}\label{th2} \noindent By the representation
(\ref{m-sacingsiden(n={k})}) we get
\begin{equation}
\left\{\gamma_{{Nk}}(t), 0\leq t <
1\right\}{\stackrel{d}{=}}\left\{\gamma_N^1(t), 0\leq t <
1\right\}.
\end{equation}
\noindent Recall that our aim is to prove the inequality
\begin{equation}
\mathbb{P}\left(\sup_{0 \leq t \leq a }|\gamma^1_{N}(t)-W^{*}_{N}(t)| >
A_{2}N^{-1/4}(\log aN)^{3/4} \right)\leq B_{2}N^{-\varepsilon},
\end{equation}
\noindent where
\begin{equation}
W^{*}_{N}(t):=B^{(1)}_N(t)-\frac{\phi_{k}(t)}{k}\int_0^{\infty}B^{(1)}_N(F_{k}(t))dt.
\end{equation}
 \noindent First, we observe that
 \begin{eqnarray}\label{inequality (4.1)}
\lefteqn{\gamma_N^1(t)-\left(B^{(1)}_N(t)-\frac{\phi_{k}
(t)}{k}\int_0^{\infty}B^{(1)}_N(F_{k}(t))dt\right)}
\nonumber\\&=&\kappa_N(t)-B^{(1)}_N(t)+\left(\left(\frac{T_N}{{Nk}}\right)^{-1}-1\right)
\kappa_N(t)\nonumber\\&&+\phi_{k}(t)N^{1/2}\left(\left(\frac{T_N}
{{Nk}}\right)-1\right)\left(\left(\frac{T_N}{{Nk}}\right)^{-1}-1\right)
\nonumber\\&&-\frac
{\phi_{k}(t)}{k}\left(N^{1/2}\left(k-\frac{T_N}{N}\right)
-\int_0^{\infty}B^{(1)}_N(F_{k}(t))dt\right).
\end{eqnarray}
\noindent Now, by Theorem \ref{1er theorem} we infer that 
\begin{equation}\label{inequality (4.2)}
\mathbb{P}\left(\sup_{0 \leq t \leq a}|\kappa_N(t)-B^{(1)}_N(t)| \geq
A_1N^{-1/4}(\log aN)^{3/4} \right)\leq B_1N^{-\varepsilon}.
\end{equation}
\noindent Noting that
\begin{equation}\label{fff1}
\sup_{0 \leq t\leq a}\phi_{k}(t)=\sup_{0\leq x\leq
Q(a)}xf_{k}(x)<\infty.
\end{equation}
\noindent Let $$A_{13}=A_{9}\sup_{0\leq x\leq Q(a)}xf_{k}(x).$$ By combining 
Lemma \ref{2eme lemme}  with (\ref{fff1}), we get
\begin{eqnarray}\label{inequality (4.3)}
\lefteqn{\mathbb{P}\left(\sup_{0\leq t\leq
a}\left|N^{1/2}\left(1-\frac{T_N}{{Nk}}\right)
\phi_{k}(t)-\frac{\phi_{k}(t)}{k}\int_0^{\infty}B^{(1)}_N\left(F_{k}(t)dt\right)\right|
> A_{13}N^{-1/2}\left(\log
aN\right)^{2}\right)}\nonumber\\&& \qquad\qquad\qquad\qquad\qquad\qquad\qquad\qquad \qquad\qquad\qquad\qquad\qquad\qquad\leq B_{9}N^{-\varepsilon}.
\end{eqnarray}
\noindent By using the following elementary fact
$$\Big(\frac{1}{u}-1\Big)=-(u-1)+\frac{1}{u}\left(u-1\right)^2, $$
we obtain  the following relation 
\begin{equation}
\left(\left(\frac{T_N}{{Nk}}\right)^{-1}-1\right)\kappa_N(t)
=-\left(\left(\frac{T_N}{{Nk}}\right)-1\right)\kappa_N(t)+
\left(\left(\frac{T_N}{{Nk}}\right)-1\right)^2\frac{{Nk}}{T_N}\kappa_N(t).
\end{equation}
\noindent We have the inequalities
\begin{eqnarray}\label{s3}
\lefteqn{\mathbb{P}\left(\sup_{0\leq t\leq
a}\left|\left(\left(\frac{T_N}{{Nk}}\right)-1\right)\kappa_N(t)\right|
>\left(A_{11}N^{-1/2}(\log aN)^{1/2}\right)\right.}\nonumber\\ &&\times \left.
\left(\left(\frac{1}{2}\varepsilon(\log aN)\right)^{1/2} +
A_{1}N^{-1/4}(\log aN)^{3/4}\right)\right)\nonumber\\ &\leq&
\mathbb{P}\left(\left|\left(\left(\frac{T_N}{{Nk}}\right)-1\right)\right|
>\left(A_{11}N^{-1/2}(\log aN)^{1/2}\right)\right)\nonumber\\&&+\mathbb{P}\left(\sup_{0\leq t\leq
a}\left|\kappa_N(t)\right|
>\left(\left(\frac{1}{2}\varepsilon(\log aN)\right)^{1/2} + A_{1}N^{-1/4}(\log
aN)^{3/4}\right)\right) \nonumber\\ &\leq&
\mathbb{P}\left(\left|\left(\left(\frac{T_N}{{Nk}}\right)-1\right)\right|
>\left(A_{11}N^{-1/2}(\log aN)^{1/2}\right)\right)\nonumber\\ &&+\mathbb{P}\left(\sup_{0\leq t\leq
a}\left|\kappa_N(t)-B_N^1(t)\right|
> \left(A_{1}N^{-1/4}(\log
aN)^{3/4}\right)\right)\nonumber\\&&+\mathbb{P}\left(\sup_{0\leq t\leq
a}\left|B^{(1)}_N(t)\right|
>\left(\frac{1}{2}\varepsilon(\log aN) \right)^{1/2}\right)\nonumber\\
&\leq& B_{11}N^{-\varepsilon}+B_1N^{-\varepsilon}+
2N^{-\varepsilon}\nonumber\\ &\leq& B_{14}N^{-\varepsilon}.
\end{eqnarray}
\noindent By the law of large numbers, $T_N/N $ converges  to ${k}$, as
$n$ tends to infinity. Then $T_N/{Nk}$ tends to one when $n$ tends
to infinity. On the other hand, we remark, if $T_N/{Nk} \geq 1/2$,
then ${Nk}/T_N \leq 2$. We can see that
\begin{eqnarray}\label{s4}
\lefteqn{\mathbb{P}\left(\sup_{0\leq t\leq
a}\left|\left(\left(\frac{T_N}{{Nk}}\right)-1\right)^2\frac{{Nk}}{T_N}\kappa_N(t)\right|
>\left(2A_{11}^2N^{-1}(\log aN)\right)\right.}\nonumber\\ &&\times \left.
\left(\left(\frac{1}{2}\varepsilon(\log aN)\right)^{1/2} +
A_{1}N^{-1/4}(\log aN)^{3/4}\right)\right)\nonumber\\ &\leq&
B_{14}N^{-\varepsilon}.
\end{eqnarray}
\noindent Using (\ref{s3}) and (\ref{s4}), we obtain
\begin{equation}\label{inequality (2.13)}
\mathbb{P}\left(\sup_{0\leq t\leq
a}\left|\left(\left(\frac{T_N}{{Nk}}\right)^{-1}-1\right)\kappa_N(t)\right|>A_{14}
N^{-1/4}(\log aN)^{3/4}\right)\leq B_{14} N^{-\varepsilon}.
\end{equation}
\noindent Moreover we have
\begin{eqnarray*}
\lefteqn
{\phi_{k}(t)N^{1/2}\left(\left(\frac{T_N}{{Nk}}\right)-1\right)\left(\left(\frac{T_N}{{Nk}}\right)^{-1}-1\right)}\\&=&
-\phi_{k}(t)N^{1/2}\left(\left(\frac{T_N}{{Nk}}\right)-1\right)^{2}
+\phi_{k}(t)N^{1/2}\left(\left(\frac{T_N}{{Nk}}\right)-1\right)^{3}\frac{T_N}{{Nk}}.
\end{eqnarray*}
\noindent Now, on $A_N(\delta)$,  we have $$\sup_{0 \leq t \leq
a}\phi_{k}(t)=M<\infty.$$Taking
$A_{15}=A^2_{11}M$ and applying similar techniques used in
line $2$ of (\ref{inequality (3.8)}) we get, by (\ref{inequality
(3.6)}), that
\begin{equation}\label{s2}
\mathbb{P}\left(\sup_{0\leq t\leq
a}\phi_{k}(t)N^{1/2}\left|\left(\left(\frac{T_N}{{Nk}}\right)-1\right)^{2}\right|>A_{15}N^{-1/2}(\log
aN) \right)\leq B_{11}N^{-\varepsilon}.
\end{equation}
\noindent Let $A_{16}=2A^3_{11}M$. By using the same
arguments, it follows 
\begin{equation}\label{s1}
\mathbb{P}\left(\sup_{0\leq t\leq
a}\phi_{k}(t)N^{1/2}\left|\left(\left(\frac{T_N}{{Nk}}\right)-1\right)^{3}\right|\frac{{Nk}}{T_N}>A_{16}
N^{-1}(\log aN)^{3/2} \right)\leq B_{11}N^{-\varepsilon}.
\end{equation}
\noindent From (\ref{s2}) and (\ref{s1}), we obtain
\begin{eqnarray}\label{inequality (2.14)}
\mathbb{P}\left(\sup_{0\leq t\leq
a}\phi_{k}(t)N^{1/2}\left(\left(\frac{T_N}{{Nk}}\right)-1\right)
\left(\left(\frac{T_N}{{Nk}}\right)^{-1}-1\right)\right.\\\left.\nonumber>
A_{17}N^{-1/2}(\log aN)\right)&\leq& B_{17}N^{-\varepsilon}.
\end{eqnarray}
\noindent Now, by combining (\ref{inequality (4.1)}),
(\ref{inequality (4.2)}), (\ref{inequality (4.3)}),
(\ref{inequality (2.13)}) and (\ref{inequality (2.14)}), we obtain 
\begin{eqnarray}\label{inequality (4.4)}
\mathbb{P}\left(\sup_{0\leq t\leq
a}\left|\gamma_N^1(t)-\left(B^{(1)}_N(t)-\frac{\phi_{k}(t)}{k}
\int_0^{\infty}B^{(1)}_N(F_{k}(t))dt\right)\right|\right.
\nonumber\\\left.>A_{2}N^{-1/4}\left(\log
aN\right)^{3/4}\right)\leq B_{2} N^{-\varepsilon}.
\end{eqnarray}
\noindent By Lemma $4.4.4$ of \cite{csorgorevesz1981} and
(\ref{gamma}), we can define a sequence of Gaussian processes
$\{W_{{Nk}}(t): 0\leq t\leq 1\}, N=1,2,\ldots$ such that for each
$N$, we have
\begin{equation}
\left\{\gamma_{{Nk}}(t),W_{{Nk}}(s): 0\leq t,s \leq
1\right\}{\stackrel{d}{=}}\left\{\gamma_N^1(t),W_N^*(t): 0\leq t
,s \leq 1\right\}.
\end{equation}
\noindent This completes the proof of Theorem \ref{2eme
theorem}.\hfill$\blacksquare$

\subsection*{Proof of Theorem \ref{3eme theorem}.}\label{th3}\vskip5pt We are going to give
the main steps of the proof. The details are the same as in
Theorem \ref{2eme theorem}. Assume first that $n={Nk}$.
Keep in mind the representation (\ref{alpha}) for the empirical process of
$k$-spacings. We are aimed   to prove the following
\begin{equation}\label{inequality (3.1)}
\mathbb{P}\left(\sup_{0 \leq x \leq Q_{k}(a)}|\alpha^1_N(x)-V^*_N(x)|
\geq A_{3}N^{-1/4}(\log aN) \right)\leq B_{3}N^{-\varepsilon},
\end{equation}
\noindent where
\begin{equation}\label{vvv*}
V^*_N(x)=B^{(1)}_N(F_{k}(x))-\frac{1}{k}xf_{k}(x)\int_0^{\infty}B^{(1)}_N(F_{k}(y))dy.
\end{equation}
\noindent By second order Taylor expansion in the
second term of (\ref{alpha}), we get
\begin{eqnarray*}
\alpha^1_N(x)-V^*_N(x)&=&\beta_N\left(x\frac{T_N}{{Nk}}\right)-B^{(1)}_{N}
\left(F_{k}\left(x\frac{T_N}{{Nk}}\right)\right)
\\&&+B^{(1)}_{N}\left(F_{k}\left(x\frac{T_N}{{Nk}}\right)\right)
-B^{(1)}_N(F_{k}(x))+N^{1/2}\left(\frac{T_N}{{Nk}}-1\right)^2x^2f'^{k}(x_N)\\
&&+\frac{xf_{k}(x)}{k}\left(N^{1/2}\left(
\frac{T_N}{{Nk}}-1\right)-\int_0^{\infty}tdB^{(1)}_N(F_{k}(t))\right),
\end{eqnarray*}
\noindent where $$|x_N-x|\leq x\left|\frac{T_N}{{Nk}}-1\right|.$$
Making use of Lemmas \ref{2eme lemme} and \ref{lemma 3}, together
with Theorem \ref{Masonvanzwet theorem} we obtain (\ref{inequality
(3.1)}). Hence together with Lemma $4.4.4$ of
\cite{csorgorevesz1981}, we can define a sequence of Gaussian
processes $\{V_{{Nk}}(x): 0\leq x< \infty\}, N=1,2,\ldots,$ such
that for each $N$ we have
\begin{equation}
\{\alpha_{{Nk}}(x), V_{{Nk}}(y): 0\leq x, y<
\infty\}{\stackrel{d}{=}}\{\alpha^1_{N}(x), V^*_{N}(y):0\leq x,
y<\infty\}.
\end{equation}
\noindent This completes the proof Theorem (\ref{3eme theorem}) for the case  where $n={Nk}$. Now, we prove the general case where $k(N-1)<n<{Nk}$.
It follows from (\ref{m-sacingsiden}) that
\begin{eqnarray}\label{inequality (3.12)}
\lefteqn{\{\alpha_{n}(x) :0\leq x<
\infty\}}\nonumber\\&{\stackrel{d}{=}}&\left\{N^{1/2}
\left(G_{N,k}\left(x\frac{S_n}{{Nk}}\right)-F_{k}(x)\right): 0\leq
x<\infty\right\},
\end{eqnarray}
\noindent where
\begin{equation}
G_{N,k}(x)=\frac{1}{N}\sum_{i=1}^{N-1}\mathds{1}_{\{Y_i<x\}}+
\frac{1}{N}\mathds{1}_{\left\{\sum_{\ell=(N-1)k+1}^{n}E_\ell<x\right\}}.
\end{equation}
\noindent Notice that we have  the following fact 
\begin{equation}\label{inequality (3.13)}
\sup_{0\leq x\leq
Q_{k}(a)}\left|G_{N,k}\left(x\frac{S_n}{{Nk}}\right)-G_{N-1}
\left(x\frac{S_n}{{Nk}}\right)\right|\leq
\frac{1}{N}+\frac{1}{N(N-1)}, 
\end{equation}
\noindent and
\begin{equation}\label{inequality (3.14)}
\mathbb{P}\left(\left|\frac{S_n}{{Nk}}-\frac{T_{N-1}}{k(N-1)}\right|>A_{18}N^{-1}(\log
aN)\right)\leq B_{18}N^{-\varepsilon}.
\end{equation}
\noindent  Set 
\begin{eqnarray}
\nonumber \mathcal{P}=\mathbb{P}\left(\sup_{0\leq x\leq
Q_{k}(a)}\left|N^{1/2}\left(G_{N,k}\left(x\frac{S_n}{{Nk}}
\right)-F_{k}(x)\right)-V_{N-1}^{*}(x)\right|\right.\\\left.
>A_{19}N^{-1/4}\left(\log
aN\right)^{3/4}\right).
\end{eqnarray}
\noindent By the use of  (\ref{inequality (3.1)}) in connection with (\ref{inequality
(3.13)}), we infer that 
\begin{eqnarray}\label{inequality (300)}
\nonumber \mathcal{P}&\leq&\mathbb{P}\left(\sup_{0\leq x\leq
Q_{k}(a)}N^{1/2}\left|F_{k}\left(x\frac{S_n}{T_{N-1}}
\right)-F_{k}(x)\right|>A_{20}N^{-1/2}(\log aN)\right)
\\\nonumber&&+\mathbb{P}\left(\sup_{0\leq x\leq
Q_{k}(a)}N^{1/2}\left|V_{N-1}^*\left(x\frac{S_n}{T_{N-1}}
\right)-V_{N-1}^{*}(x)\right|>A_{21}N^{-1/2}(\log
aN)\right)\\&&+B_{3} N^{-\varepsilon}.
\end{eqnarray}
\noindent Once more, by a first order the Taylor expansion, we have 
\begin{equation}
N^{1/2}\left|F_{k}\left(x\frac{S_n}{T_{N-1}}\right)-F_{k}(x)
\right|=xf_{k}(x_N).N^{1/2}\left|\frac{S_n-T_{N-1}}{T_{N-1}}\right|,
\end{equation}
\noindent where $$|x_N-x|\leq
x\left|\frac{S_n-T_{N-1}}{T_{N-1}}\right|.$$ 
By combining 
Lemma \ref{2eme lemme} with  (\ref{inequality (3.14)}), it follows that
\begin{equation}\label{inequality (3.15)}
\mathbb{P}\left(\left|\frac{S_n}{T_{N-1}}-1\right|>A_{22}N^{-1}(\log
aN)\right)\leq B_{22}N^{-\varepsilon}.
\end{equation}
\noindent By arguing in a similar way as in the proof
(\ref{inequality (3.8)}), we obtain that
\begin{equation}\label{inequality (3.16)}
\mathbb{P}\left(\sup_{0\leq x\leq
Q_{k}(a)}N^{1/2}\left|F_{k}\left(x\frac{S_n}{T_{N-1}}\right)-F_{k}(x)\right|>A_{20}N^{-1/2}(\log
aN)\right)\leq B_{20}N^{-\varepsilon}.
\end{equation}
\noindent Now, by definitions (\ref{vvv*}), (\ref{inequality
(3.16)}), and through a similar argument as that used in the end
of the proof of Lemma \ref{lemma 3}, we get
\begin{equation}\label{inequality (A39)}
\mathbb{P}\left(\sup_{0\leq x\leq
Q_{k}(a)}N^{1/2}\left|V_{N-1}^*\left(x\frac{S_n}{T_{N-1}}\right)-V_{N-1}^{*}(x)\right|>A_{21}N^{-1/2}(\log
aN)\right)\leq B_{21}N^{-\varepsilon}.
\end{equation}
\noindent Then, making use of the equations  (\ref{inequality (300)}),(\ref{inequality (3.16)}) and  (\ref{inequality
(A39)}),   we obtain 
\begin{eqnarray}
\mathbb{P}\left(\sup_{0\leq x\leq
Q_{k}(a)}\left|N^{1/2}\left(G_{N,k}\left(x\frac{S_n}{{Nk}}\right)-F_{k}(x)\right)-V_{N-1}^{*}(x)\right|
\right.\\\left.>A_{23}N^{-1/4}\left(\log
aN\right)^{3/4}\right)\leq B_{23} N^{-\varepsilon}.
\end{eqnarray}
\noindent Again, by Lemma $4.4.4$ of \cite{csorgorevesz1981} and
(\ref{inequality (3.12)}), we can get a sequence of Gaussian
processes $\{V_n(x): 0\leq x<\infty\}, k(N-1)<n<{Nk}, N=1,2\ldots$,
such that for each $N$ we have
\begin{eqnarray*}
\lefteqn{\{\alpha_{n}(x), V_{n}(y): 0\leq x, y<
\infty\}}\\&{\stackrel{d}{=}}&\left\{N^{1/2}\left(G_{N,k}\left(x\frac{S_n}{{Nk}}\right)-F_{k}(x)\right),
V^*_{N-1}(y): 0\leq x, y<\infty\right\}.
\end{eqnarray*}
This completes the proof of Theorem \ref{3eme
theorem}.\hfill$\blacksquare$\vskip5pt
\subsection*{Proof of Corollary~\ref{corol22}}

The functional $\Phi$ being Lipschitz, there exists a positive constant $L$
such that, for any functions $v,w$,
\begin{equation}\label{lip}
|\Phi(v)-\Phi(w)|\leq L \sup_{t\in\mathbb{R}} |v(t)-w(t)|.
\end{equation}
Let us choose for $v,w$ the processes $$\gamma_n:=\gamma_{{Nk}}(\cdot) \quad\mbox{and}\quad
\mathbb{W}_n:=W_{{Nk}}(\cdot).$$
Applying the elementary inequality
$$
|\mathbb{P}(A)-\mathbb{P}(B)|
\leq \mathbb{P}(A\backslash B) + \mathbb{P}(B\backslash A)
=\mathbb{P}\big((A\backslash B)\cup (B\backslash A)\big)
$$
to the events $A=\{\Phi(\gamma_n)\leq x\}$ and $B=\{\Phi(\mathbb{W}_n)\leq x\}$ provides,
for any $x\in\mathbb{R}$ and any $n\in\mathbb{N}^*$,
$$
\big|\mathbb{P} \{\Phi(\gamma_n)\leq x\}-\mathbb{P} \{\Phi(\mathbb{W}_n)\leq x\}\big|
\leq \mathbb{P} \{\Phi(\gamma_n)\leq x\leq \Phi(\mathbb{W}_n) \quad\mbox{or}\quad\Phi(\mathbb{W}_n)\leq x\leq \Phi(V_n)\}.
$$
Now, applying to the elementary fact that [$a\le x\le b$ or $b\le x\le a$
implies $|b-x|\le |b-a|$] to the numbers $a=\Phi(\gamma_n)$ and $b=\Phi(\mathbb{W}_n)$
\begin{align*}
\mathbb{P} \{\Phi(\gamma_n)\leq x\leq \Phi(\mathbb{W}_n)\}&
\leq\mathbb{P} \!\left\{|\Phi(\mathbb{W}_n)-x|\leq |\Phi(\mathbb{W}_n)-\Phi(\gamma_n)|\right\}\!,
\end{align*}
from which, due to (\ref{lip}), we deduce that
\begin{equation}\label{estimPhiVW}
\big|\mathbb{P} \{\Phi(\gamma_n)\leq x\}-\mathbb{P} \{\Phi(\mathbb{W}_n)\leq x\}\big|
\leq \mathbb{P} \!\left\{|\Phi(\mathbb{W}_n)-x| \leq L \sup_{t\in\mathbb{R}} |\gamma_n(t)-\mathbb{W}_n(t)|\right\}\!.
\end{equation}
We obtain the estimate below valid for large enough $n$, for $\epsilon_n=c(N^{-1/4}(\log aN)^{3/4})$
\begin{equation}\label{estimationbis}
\mathbb{P}\!\left\{ \sup_{t\in\mathbb{R}} \left|\gamma_n(t)-\mathbb{W}_n(t)\right|
\geq \epsilon_n\right\}\leq 
o\left(N^{-1/4}(\log aN)^{3/4}\right)
\end{equation}
Now, by (\ref{estimPhiVW}), we write
\begin{eqnarray}
\nonumber\lefteqn{
\big|\mathbb{P} \{\Phi(\gamma_n)\leq x\}-\mathbb{P} \{\Phi(\mathbb{W}_n)\leq x\}\big|}\\
&\leq&\nonumber
\mathbb{P} \!\left\{ \sup_{t\in\mathbb{R}}
\left|\gamma_n(t)-\mathbb{W}_n(t)\right|<\epsilon_n, |\Phi(\mathbb{W}_n)-x| \leq
L \sup_{t\in\mathbb{R}} |\gamma_n(t)-\mathbb{W}_n(t)|\right\}
\nonumber\\
\nonumber&
&+\mathbb{P}\!\left\{ \sup_{t\in\mathbb{R}} \left|\gamma_n(t)-\mathbb{W}_n(t)\right|
\geq \epsilon_n,|\Phi(\mathbb{W}_n)-x| \leq
L \sup_{t\in\mathbb{R}} |\gamma_n(t)-\mathbb{W}_n(t)|\right\}
\nonumber\\
&\leq&\;
\mathbb{P} \!\left\{|\Phi(\mathbb{W}_n)-x| \leq L \epsilon_n\right\}
+\mathbb{P}\!\left\{ \sup_{t\in\mathbb{R}} \left|\gamma_n(t)-\mathbb{W}_n(t)\right|
\geq \epsilon_n\right\}\!.
\label{estimationter}
\end{eqnarray}
Noticing that the distribution of $W_{{Nk}}$ does not depend on $n$,
which entails the equality
$$
\mathbb{P} \!\left\{|\Phi(\mathbb{W}_n)-x| \leq L \epsilon_n\right\}
=\mathbb{P} \!\left\{|\Phi(\mathbb{W})-x| \leq L \epsilon_n\right\}
$$
where $\mathbb{W}$ is a gaussian process with the same distribution as $W_{{Nk}}$.
and recalling the assumption that the r.v.
$\Phi(\mathbb{W})$ admits a density function bounded by $M$ say,
we get that, for any $x\in\mathbb{R}$ and any $n\in\mathbb{N}^*$,
\begin{equation}\label{estimationquater}
\mathbb{P} \!\left\{|\Phi(\mathbb{W}_n)-x| \leq L \epsilon_n\right\}\leq 2LM\epsilon_n.
\end{equation}
Finally, putting (\ref{estimationbis}) and  (\ref{estimationquater})
into (\ref{estimationter}) leads to (\ref{estimationfunctional}),
which completes the proof of Corollary~\ref{corol22}.
An alternative proof of a similar result may be found in \cite{ShorackWellner1986} pp. 502--503.
\hfill$\blacksquare$

\subsection*{Proof of the Theorem \ref{theoreme increment m spacings}}

\noindent The proof of this result is based on the Theorem \ref{theoreme1.2}
and the Lemma \ref{1er lemmme chap2} below. Let us consider the restriction
on the interval $[0,1]$ of a standard Wiener process $\{W(t):t\geq 0\}$.

\begin{lemma}
\label{1er lemmme chap2} Assume that $\{h_{N}:N\geq 1\}$ verify \textbf{(H.1)} and 
\textbf{(H.2)}. Then 
\begin{equation}
\lim_{N\rightarrow \infty }\Delta \left( \{\beta
_{N}(W(t+h_{N}\cdot)-W(t)):0\leq t\leq 1-h_{N}\},\mathcal{S}_{\frac{c}{c+1}
}\right) \overset{\mathbb{P}}{=}0.  \label{eq1}
\end{equation}
\end{lemma}

\noindent \textbf{Proof of Lemma \ref{1er lemmme chap2} .} For $c<\infty $, the proof of this lemma
is due to the joint use of a result of \cite
{DeheuvelsRevesz1993} and the scaling property. Notice that the result of  \cite
{DeheuvelsRevesz1993} is obtained with our the hypothesis of monotonically
decreasing (resp. increasing) when $h_{N}\rightarrow 0$ (resp. $
Nh_{N}\rightarrow \infty $), we may refer also to  \cite{DeheuvelsRevesz1986}. For the case of a non bounded arbitrary $c$, the proof of  \cite{DeheuvelsRevesz1993} could still be used. We just notice that the quantity $\delta _{T}$ defined in their article
tends again to $0$. Then, the replacement of their Theorem $1.1$ initially
used, by the Theorem $1.2$ appearing in their article in connection with the use of the
property of scaling leads to the expected result.
\hfill$\blacksquare$
\vskip5pt \noindent Now, we are equipped  to prove
Theorem \ref{theoreme increment m spacings}. \subsection*{Proof of Theorem \ref{theoreme increment m spacings}.} For all $
0\leq t\leq 1-h_{N}$, we define the quantity, 
\begin{equation}
\omega _{N}(h_{N},t;\cdot)=W(t+h_{N}\cdot)-W(t).
\end{equation}
Through (\ref{eq1}), the triangle inequality applied on the Hausdorff
distance bring our proof to show that, as $n\rightarrow \infty $, 
\begin{equation*}
\boldsymbol{\Delta} \Big(\Big\{\beta _{N}\eta _{N}(h_{N},t;\cdot):0\leq t\leq 1-h_{N}\Big\},
\Big\{\beta _{N}\omega _{N}(h_{N},t;\cdot):0\leq t\leq \Gamma _{1}^{\prime
}-h_{N}\Big\}\Big)\overset{\mathbb{P}}{\rightarrow }0,
\end{equation*}
By (\ref{eq0}) and (\ref
{eq2.1 AlyBeirlHor84}), it is easy to remark  that, for all $\varepsilon >0$,
there exists a constant $C_{\varepsilon }$ such that the following
inequality takes place for all $n$ large enough, with probability greater 
than $1-\varepsilon $, 
\begin{eqnarray}
\lefteqn{\Big\|\beta _{N}\eta _{N}(h_{N},t;\cdot)-\beta _{N}\omega _{N}(h_{N},t;\cdot)\Big\|}\nonumber \\
&\leq& \beta _{N}C_{\varepsilon }n^{-1/4}(\log n)^{3/4}  \label{eq3} 
+\Big\|\beta _{N}(\Gamma _{N}^{\prime }(t+h_{N}\cdot)-\Gamma _{N}^{\prime
}(t))-\beta _{N}\omega _{N}(h_{N},t;\cdot)\Big\|.
\end{eqnarray}
The hypotheses \textbf{(H.{2})} and \textbf{(H.{3})} implies that, as $n\rightarrow \infty $, 
\begin{equation*}
\beta _{N}n^{-1/4}(\log n)^{3/4}\rightarrow 0,
\end{equation*}
Let us consider now the 
Brownian bridge 
\begin{equation*}
B(u)=W(u)-uW(1),0\leq u\leq 1.
\end{equation*}
Recall  the definition (\ref{eq2.1 AlyBeirlHor84}) of $\Gamma _{N}^{\prime }$,  we infer that we have 
\begin{eqnarray}\nonumber
\lefteqn{\Big\|\beta _{N}(\Gamma _{N}^{\prime }(t+h_{N}\cdot)-\Gamma _{N}^{\prime
}(t))-\beta _{N}\omega _{N}(h_{N},t;\cdot)\Big\|}
\\
&\leq& \beta _{N}h_{N}|W(1)|
\label{eq4} 
+\beta _{N}\Big\|\{\phi _{k}(t)-\phi
_{k}(t+h_{N}\cdot)\}\int_{0}^{1}B(u)dQ_{k}(u)\Big\|,
\end{eqnarray}
where we recall $$\phi _{k}(t)=\frac{1}{k}Q_{k}(t)f_{k}(Q_{k}(t)),$$\vskip5pt\noindent
and $$\sup_{0\leq t\leq 1}\phi _{k}(t)<\infty ,$$ which implies that 
\begin{equation}
\beta _{N}\Big\|\{\phi _{k}(t)-\phi _{k}(t+h_{N}\cdot)\}\Big\|=o(1)\times
O_{\mathbb{P}}(1)\overset{\mathbb{P}}{\rightarrow }0~~\mbox{as}~~n\rightarrow \infty .
\label{eq5}
\end{equation}
It follows from \textbf{(H.1)} that $\beta _{N}h_{N}\rightarrow 0$ as $
n\rightarrow \infty $. 
By  (\ref{eq4}) and (\ref{eq5}) and making use of Tchebychev's inequality, we obtain readily that, as $n\rightarrow \infty $, 
\begin{equation}
\Big\|\beta _{N}\eta _{N}(h_{N},t;\cdot)-\beta _{N}\omega _{N}(h_{N},t;\cdot)\Big\|
\overset{\mathbb{P}}{\rightarrow }0,
\end{equation}
Hence, we  finally obtain, as $n\rightarrow \infty $
\begin{equation*}
\Delta \Big(\Big\{\beta _{N}\eta _{N}(h_{N},t;\cdot):0\leq t\leq 1-h_{N}\Big\},
\Big\{\beta _{N}\omega _{N}(h_{N},t;\cdot):0\leq t\leq 1-h_{N}\Big\}\Big)\overset
{\mathbb{P}}{\rightarrow }0,
\end{equation*}
The proof of Theorem \ref{theoreme increment m spacings} is therefore completed.
\hfill$\blacksquare$

\subsection*{Proof of Corollary \ref{corollaire1}.} As $\mathcal{S}_{c/(c+1)}$ is a compact
and connected set of $B_{0}[0,1]$ equipped with the uniform topology,
the continuity of the functional $\Phi $ implies that image set $\Phi (\mathcal{S}
_{c/(c+1)})$  is also a compact and connected set of $\mathbb{R}$ of the form $[l_{\Phi },L_{\Phi }]$, where $$l_{\Phi }=\inf_{f\in \mathcal{S
}_{c/(c+1)}}\Phi (f) ~~\mbox{and}~~ L_{\Phi }=\sup_{f\in \mathcal{S}_{c/(c+1)}}\Phi
(f).$$ By making use of Theorem \ref{theoreme increment m spacings} and the continuity of the functional $
\Phi $,  we obtain 
\begin{equation*}
\lim_{n\rightarrow \infty }\boldsymbol{\Delta} \Big(\big\{\beta _{N}\eta
_{N}(h_{N},t;\cdot):0\leq t\leq 1-h_{N}\big\},[l_{\Phi },L_{\Phi }]\Big)\overset{
\mathbb{P}}{=}0.
\end{equation*}
\noindent which readily imply  that
\begin{equation*}
\lim_{n\rightarrow \infty }\sup_{0\leq t\leq 1-h_{N}}\Phi (\beta _{N}\eta
_{N}(h_{N},t;\cdot))\overset{\mathbb{P}}{=}L_{\Phi }.
\end{equation*}
Hence the proof of Corollary \ref{corollaire1} is complete. \hfill$\blacksquare$\subsection*{Proof of Corollary \ref{corollaire2}.} For any $f(\cdot)\in \mathcal{S}_{c/(c+1)}$,  the Schwarz inequality shows that, for each $0 \leq s \leq 1$,
\begin{equation*}
|f(s)|=\left\vert \int_{0}^{s}\dot{f}(u)du\right\vert \leq \left(
s\int_{0}^{1}\dot{f}^{2}(u)du\right) ^{1/2}\leq \left( \frac{c}{c+1}\right)
^{1/2}.
\end{equation*}
Therefore, by choosing $f(\cdot)$ in such a way that 
\begin{equation*}
f(s)=\left( \frac{c}{c+1}\right)^{1/2}s,~\mbox{for}~0\leq s\leq 1,
\end{equation*}
gives, after a simple calculation,  that $$f(\cdot)\in \mathcal{S}_{c/(c+1)}~\mbox{ and }~ \Phi
_{1}(f)=\Phi _{2}(f)=\left( \frac{c}{c+1}\right)^{1/2}.$$Thus, 
 an application of Corollary \ref{corollaire1} implies, in turn, that,  (\ref{eq1 corrolaire2}) and (\ref{eq2
corrolaire2}) hold. Moreover, for all  $f(\cdot)\in \mathcal{S}
_{c/(c+1)}$, using Schwarz' inequality implies readily that  
\begin{equation*}
\mp \int_{0}^{1}v(u)\dot{f}(u)du\leq \left( \frac{c}{c+1}
\int_{0}^{1}v^{2}(u)du\right) ^{1/2}.
\end{equation*}
The function $f(\cdot)$ of derivative, 
\begin{equation*}
\dot{f}(s)=\mp \left( \frac{c}{c+1}\right)^{1/2}\left( \int_{0}^{1}v^{2}(u)du\right)
^{-1/2}v(s),
\end{equation*}
belongs to $\mathcal{S}_{c/(c+1)}$. Furthermore, 
\begin{equation*}
\Phi _{3}(f)=\left( \frac{c}{c+1}\right)^{1/2}\left( \int_{0}^{1}v^{2}(s)ds\right) ^{1/2}.
\end{equation*}
An application of Corollary \ref{corollaire1} implies that  (\ref{eq3 corrolaire2}) hods..\hfill$\blacksquare$

\end{document}